\documentclass[10pt]{article}
\usepackage{layout}
\usepackage{amssymb}
\usepackage{latexsym}
\usepackage{amsmath}
\usepackage{color}
\usepackage[T1]{fontenc}

\newcommand{\Ref}[1] {(\ref{#1})}
\newcommand{\adelta}{(-\Delta)^{\alpha}}
\newcommand{\R}{\mathbb{R}}

\newcommand{\N}{\mathbb{N}}
\newcommand{\T}{\mathbb{T}^{2}}
\newcommand{\Z}{\mathbb{Z}}
\newcommand{\be}{\begin{equation}}
\newcommand{\ee}{\end{equation}}
\newcommand{\fin}{\rule{3mm}{0mm} \hfill \rule{2mm}{2mm}}
\newcommand{\non}{\nonumber}

\newtheorem{thm}{Theorem}

\newtheorem{lem}{Lemma}
\newtheorem{prop}{Proposition}

\newtheorem{rem}{Remark}

\newcommand{\edoc}      {\end{document}}

\newcommand{\bean}     {\begin{eqnarray}\nonumber}
\newcommand{\bea}     {\begin{eqnarray}}
\newcommand{\blll}     {\begin{array}{lll}}
\newcommand{\brcl}     {\begin{array}{rcl}}
\newcommand{\barr}     {\begin{array}}
\newcommand{\earr}     {\end{array}}
\newcommand{\eea}     {\end{eqnarray}}
\newcommand{\beet}     {\begin{eqnarray*}}
\newcommand{\eec}     {\end{center}}
\newcommand{\bequ}      {\begin{equation}}
\newcommand{\eequ}      {\end{equation}}
\newcommand{\bet}     {\begin{tabular}}
\newcommand{\eet}     {\end{tabular}}
\newcommand{\btab}     {\begin{table}}
\newcommand{\etab}     {\end{table}}
\newcommand{\ben}     {\begin{enumerate}}
\newcommand{\een}     {\end{enumerate}}
\newcommand{\bec}     {\begin{center}}
\newcommand{\beit}     {\begin{itemize}}
\newcommand{\eeit}     {\end{itemize}}
\begin{document}


\title{\scshape A SEMI DISCRETE DYNAMICAL SYSTEM FOR A 2D DISSIPATIVE QUASI GEOSTROPHIC EQUATION }
\author{M. MOALLA-TRABELSI  and  E. ZAHROUNI  \\
 Unit\'e de recherche : Multi-Fractals et Ondelettes  \\
 Facult\'e des Sciences de Monastir\\
Av. de l'environnemnt, Monastir\\
 Tunisia} 

\maketitle

\bigskip

\begin{abstract}
A semi-discretization in time, according to a full implicit Euler
scheme, for a 2D dissipative quasi geostrophic equation, is studied.
We prove existence, uniqueness and regularity results of the
solution to the predicted discretization, in the subcritical case
for any initial data in $\dot{L}^2$. Hence, we define an infinite
semi-discrete dynamical system, then we prove the existence and the
regularity of the corresponding global attractor, for a source term
$f$ in $\dot{L}^{p_{\alpha}}$, for a fixed $p_{\alpha}  = \frac{2}{ 1-\alpha} $.
\end{abstract}

\smallskip

\section{Introduction}

\noindent In this paper, we focus on a two dimensional dissipative
quasi-geostrophic equation (QG),
\begin{eqnarray}
\label{qg1}
\partial_t{\theta} + \nu ( - \Delta)^{\alpha}\theta  +  u.\nabla \theta   =
f.
\end{eqnarray}

\noindent The solution $\theta $ of \Ref{qg1} is a real valued
function defined on $\R_+ \times \Omega, $ where $\Omega$ is either
$\R^2$ or $\T = ] 0 , 2 \pi [^2.$ We assume that $\theta$ satisfies
the following initial condition:
\begin{eqnarray}
\label{qg1a} \theta(0,x)  =  \theta_{0}(x).
\end{eqnarray}

\noindent The solution $\theta$ represents the temperature of the
fluid and $u=(u_1,u_2)$ is the divergence free velocity field which
is related to $\theta$ by the mean of Riesz transforms according to:
\bea\label{velocity} u \; = \;{\mathcal R^{\bot}}(\theta) \;
 = \; (-\Lambda^{-1}{\mathcal R}_2\theta, \Lambda^{-1}{\mathcal R}_1\theta )
= ( -\partial_{x_2} (-\Delta)^{-\frac 12} \theta, \partial_{x_1}
(-\Delta)^{-\frac 12} \theta
 ).\eea

\noindent The source term $f$ is at least square integrable and time
independent. $ \nu > 0 $ is the viscosity coefficient, and $\alpha
\in ( 0 , 1 ) $ is a fixed parameter. In the case where $\Omega =
\T$, we suppose that \bea \non \theta \quad \hbox{is} \; 2\pi \;
\hbox{periodic}\; \hbox{in each direction.} \eea

We notice that the case $\alpha=\frac 12$ is the dimensionally
correct analogue of the 3D Navier-Stokes equation, this case is
therefore called the critical one. Then $\alpha>\frac 12$ is called
the subcritical case and $0<\alpha<\frac 12$ is the supercritical
one. These models arise under the assumptions of fast rotation,
uniform stratification and uniform potential vorticity. The reader
is referred to Canstantin, Majda and Tabak \cite{consmajtab}, Held
and {\it al} \cite{Heldal}, Pedolsky \cite{pedolsky} and the
references therein for more details.

Nowadays, there are intensive investigation about existence,
uniqueness and regularity of solutions of \Ref{qg1} in the
continuous case, for different values of the diffusion parameter
$\alpha$. Indeed, since the pioneering work of S. Resnick in
\cite{resnick}, where weak solutions have been constructed, we can
cite the works of Constantin C\'ordoba and Wu \cite{const2001} in
the critical case, and the one of Chae and Lee \cite{chae2003}, in
both critical and supercritical cases.

\noindent Let us mention also the work of Kiselev, Nazarov, and
Volberg \cite{kiselev}, and the work of Cafarelli and Vasseur
\cite{caffarelli}, in the same direction.

\noindent We focus here on the subcritical case. In this framework,
we refer the reader essentially to the paper of Constantin and Wu
\cite{constw}, where the authors showed that any solution with
smooth initial value is smooth for all time. On the other hand, long
time behavior of solutions of \Ref{qg1}, was studied by N. Ju in
\cite{Ju2} and by Berselli in \cite{Berselli}. These results give
the proof of the existence of a global attractor, for the semi-group
generated by the solutions of the quasi-geostrophic equations. An
interesting question, is whether or not this important dynamical
behavior, can be captured by some of the classical numerical schemes
for solving (QG).

In fact, from the numerical point of view, Constantin and {\it al}
in \cite{constal}, performed a careful numerical study of the long
time behavior of solutions to (QG).

To the best of our Knowledge, numerical schemes for solving
\Ref{qg1} are seldom studied or even non-existent. Then, we focus in
this paper, on a semi-discretization in time of \Ref{qg1} according
to a full implicit Euler scheme, keeping the space variable
continuous. In particular, we opted for this commonly used Euler
scheme, since it is known to be a first order convergence scheme,
unconditionally stable. Then, this scheme seems to be a suitable and
reliable one to give answer to our expectations.

Here we are concerned with the discrete dynamical system associated
to \Ref{imp1}. More precisely, we prove firstly existence uniqueness
and regularity of the solution to \Ref{imp1}, and secondly we prove
the existence of the global attractor. 


\vspace{0.5cm} We notice that N. Ju considered in \cite{Ju3}, a time
discretization of the non-stationary viscous incompressible
 Navier-Stokes equations, according to a linear backward Euler scheme.
 He treated either the no-slip boundary condition or the periodic boundary
condition, in a 2D bounded domain, with a non-zero external force.
 As one of the main results obtained in \cite{Ju3},
 the global attractor for the approximation' scheme was proved to
exist.



\vspace{0.5cm}

\noindent This paper is organized as follows. In section 2, we set
our framework and state the main results.

\noindent Section 3 is devoted to prove the Theorem \ref{TH1} which
states existence, uniqueness and regularity results of a solution to
the discretized scheme. 

\noindent Finally, in section 4, we prove the Theorem \ref{TH2},
namely the existence and the regularity of the global attractor.

\section{The framework and main results}

In this section, we review the notations used throughout the
article, and we refer to some mathematical tools, which are very
useful to the success of our discussion. We set $\Omega=\T$, and let
$L^p(\Omega)$ denotes the space of the pth-power integrable
functions normed by
$$ \parallel f  \parallel_ p =(\int |f (x)|^p dx)^{\frac 1p}$$
for $p \in [1,\infty).\;$ \noindent As usual, $\hat{f}$ is the
Fourier transform of f, i.e.
 \bea \non\hat{f}(k) = \frac{1}{(2\pi)^2}\int_{\Omega} f(x) e^{-ik.x} dx.
\eea $\Lambda= (-\Delta)^{\frac 12},$ denotes the
pseudo-differential operator given by \bea
 \non
 \hat{(\Lambda f)}(k) =
|k| \hat{f}(k). \eea More generally, $\Lambda f$  can be identified
by means of Fourier series as, \bea
 \non
\Lambda^{\beta}f(x) = \sum_{k \in \Z^2} |k|^{\beta} \hat{f}(k)
e^{ik.x}.
 \eea
We define the Sobolev spaces
 \bea
  \non
  H^{s,p}= H^{s,p}(\Omega) = \left\{f \in L^p(\Omega),\quad \Lambda^s f \in L^p(\Omega)\right\} .
  \eea
\noindent Since we consider periodic boundary conditions on
$\Omega$, obviously, all derivatives of the solution $\theta$ are
mean zero. Then, $\bar{\theta}$ the mean value of $\theta$ satisfies
$$ \frac{d}{dt}\bar{\theta} = \frac{1}{\mid\Omega\mid} \frac{d}{dt}\int_\Omega \theta dx =\bar{f},$$
where
$$\bar{\theta}= \frac{1}{\mid\Omega\mid}\int_\Omega \theta(x) dx,\;\;and\;\;\;\bar{f}= \frac{1}{\mid\Omega\mid}\int_\Omega f(x)
dx.$$ Hence, without loss of generality, we may restrict the
discussion to $\theta$ that obeys for all time to $\bar{\theta} =
0$. Otherwise, we can replace $f$ with $f-\bar{f}$ and $\theta$ with
$\theta-\bar{\theta}$ and Eq.\Ref{qg1} will not change essentially.
For that purpose we introduce
 \bea  \non \dot{H}^{s,p}=
\dot{H}^{s,p}(\Omega) = \left\{f \in H^{s,p},\quad \int_ {\Omega}
f(x)  dx = 0 \right\}. \eea Accordingly, we introduce the sobolev
spaces $$ \dot{H}^{s}= \dot{H}^{s,2}\qquad \hbox{and} \quad
\dot{L}^p= \dot{H}^{0,p}$$

\noindent In order to consider (QG) as a dynamical system, we assume
that f is a given time independent scalar function, which belongs at
least to $\dot{L}^2$. We assume that the initial data $\theta_0$
belongs to $L^2$ and satisfies, \bea
\label{nulle}\frac{1}{\mid\Omega\mid} \int_{\Omega} \theta_0(x) dx =
0, \eea
so that the solution $\theta$ satisfies also \Ref{nulle}.\\

\noindent Now, we are in position to introduce the numerical scheme.
\noindent Let $\tau>0$ be a fixed real, and set $t^n=n\tau$ for
$n\in\N$. Now we recursively construct elements $\theta^{n+1}$
which approches $\theta(t^{n+1})$, by setting:  
%
%

%
%
\bea \label{imp1} \frac{\theta^{n+1}-\theta^n}{\tau} + \nu
(-\Delta)^{\alpha} \theta^{n+1} +\nabla.( u^{n+1}\theta^{n+1}) = f.
\eea

\noindent Notice that $\theta^{0}$ is an approximation of
$\theta_{0}$, and $u^{n+1} \; =  \; {\mathcal
R^{\bot}}(\theta^{n+1})$.

\noindent In our study, we consider the space domain $\Omega = \T$
and we follow the guidelines of \cite{Ju3}. Moreover, we take use of
the strategy of N. Ju in \cite{Ju2}, we quote particularly the
improved positivity lemma proved in \cite{Ju2}, which is of major
utility for our success at this stage, supplemented by the
generalized commutator estimate due to Kenig Ponce and Vega
\cite{kpv}.

\noindent Our main results state as follows :

\begin{thm}\label{TH1}
Let $\alpha \in ] 0 , 1 ]\; $ and $\; f \in  \dot{L}^2.$ Then, for
all $\theta^n \in \dot{L}^2\;$ there exists at least one solution
$\theta^{n+1}$ of \Ref{imp1} which belongs to $\dot{H}^{\alpha}$.
Moreover, if $\frac 23<\alpha < 1$ then $ \theta^{n+1}  \in
\dot{H}^{2\alpha}$ and when $\tau > 0$ is small enough, and $\frac
23\leq\alpha < 1$   this solution is unique.
\end{thm}
\noindent Furthermore, let \bea H :=  \left\{\theta\in
\dot{L}^{p_\alpha};\;
\parallel \theta
\parallel_{p_{\alpha}}  \leq M  \right\},\eea
where $M > 0$ is conveniently chosen, and consider the map $$S : H
\rightarrow H,\;\quad \theta^{n} \mapsto \theta^{n+1},$$ defined by
\Ref{imp1}. We denote by $d$ the metric distance defined by the
$\dot{L}^2$ norm, then we state:
\begin{thm}\label{TH2}
Let $\frac 23 <\alpha< 1$  and suppose that $ \;f\; \in
\dot{L}^{p_{\alpha}}$ with $p_{\alpha}=\frac{2}{1-\alpha}$. Then the
map $\; S : H \rightarrow H \; $ is continuous with respect to the
$\dot{L}^2$ topology and defines a discrete dynamical system
$(S^n)_n$ on the  complete metric space $(H,\;d)$. Besides,
$(S^n)_n$  possesses a global attractor $\mathcal{A}$ in $H,$ which
is a compact subset in $\dot{H}^{\alpha}$ and included in
$\;\dot{H}^{2\alpha}.$
\end{thm}

\vspace{0.5cm}\noindent Actually, in order to prove the above
results, we enounce different lemmas and inequalities used in the
later proofs. Let us start by a technical Lemma, which is a
consequence of the Brouwer's lemma \cite[ p.164]{temam}:
\begin{lem}\label{BrLem}
Let $X$ be a finite dimensional Hilbert space endowed with the inner
product $( .\; , \; .)$ and with the corresponding norm $
\parallel \;. \;\parallel$, and set $F:\; X \rightarrow X$
a continuous form that satisfies: $ \exists\; R>0\; $ such that,
$$[ F(\xi) , \xi ] \geq 0 \; \;\;\hbox{for}\;\; |\xi| \leq R, $$
then, there exists $\xi_0 \in X$ such that $\; |\xi_0| \leq R,\;$  and $ F(\xi_0) = 0. $\\
\fin \\
\end{lem}

\noindent Next, let us refer to the work of N. Ju, to recall an improved positivity Lemma \cite[Lemma 3.3 page 167]{Ju2},\\


\noindent\begin{lem}\label{lemlp} Let $ p \geq 2,\; s \in [ 0, 2 ],$
and $\Omega=\mathbb{T}^2$, then suppose that $\theta$ the solution
of \Ref{qg1} belongs to $L^p$, and so is $\Lambda^s \theta$. Then we
have \bea \label{positivity} p \int_{\Omega} |\theta|^{p-2} \theta
\Lambda^s \theta \geq 2 \int_{\Omega} (\Lambda^{\frac s2}
|\theta|^{\frac p2}) ^2. \eea
\end{lem}

\fin \\

\noindent The uniform Gronwall lemma presented in Temam
\cite{temam2}, is a powerful tool for a priori estimation. We recall
a discrete version of the uniform Gronwall lemmas given in Shen
\cite{shen}, which will be useful in our discussion.

\noindent\begin{lem}\label{DUGL}({\bf Discrete Uniform Gronwall Lemma.}) \\
Let $\Delta t>0$ and let $(f_n)$,$(g_n)$ and $(y_n)$ be three
positive sequences. Suppose that $\exists n_0 \geq 0,\quad r>0,\quad
a_0(r),\quad a_1(r),\quad a_2(r)$ non negative functions such that
\bea\nonumber \frac{y_{n+1}-y_{n}}{\Delta t}\leq f_n y_n + g_n,
\qquad \forall n\geq 0,\;\;\;\;\;\;\;\;\;\;\;\;\;\;\;\;\;\;\;\;\;\;\;\;\;\;\;\;\;\;\;\;\;\\
\nonumber \forall k_0\geq n_0 \qquad \Delta t\sum_{n=k_0}^{N+k_0}
f_n\leq a_0(r) ;\quad \Delta t\sum_{n=k_0}^{N+k_0} g_n\leq a_1(r);
\quad \Delta t\sum_{n=k_0}^{N+k_0} y_n\leq a_2(r),\eea $N=
[\frac{r}{\Delta t}]$ an integer, then \bea \nonumber y_n\leq (a_1 +
\frac{a_2}{r}) \exp (a_0) \qquad \forall n\geq n_0 + N. \eea
\end{lem}

\fin\\

%
%

\noindent Let us also recall a product estimate in Sobolev spaces,
due to Kenig, Ponce and Vega \cite{kpv},
\begin{equation}\label{KPV}
\parallel\Lambda^s (u\theta)\parallel_2\leq
c[\parallel u\parallel_q\parallel\Lambda^s \theta\parallel_p +
\parallel \theta\parallel_q \parallel\Lambda^s
u\parallel_p],
\end{equation}
where $\frac 1p + \frac 1q = \frac 12, $ and  the Poincaré's
inequality,  \bea\label{poincaré}
\parallel \Lambda^{\alpha}\theta\parallel_{2}^{2}&\geq
& C_0 \parallel \theta\parallel_{2}^{2},\eea where denoting by
$\lambda_1$ the first nonnegative eigenvalue of the operator
$(-\Delta)$, with periodic boundary conditions, then $C_0=
\lambda_{1}^{\alpha}$.\\
\noindent Finally, we introduce the Faedo-Galerkin method used for
the resolution of nonlinear variational formulations. That is, for
any $ \; m \in \N^{*},\;$ we consider the finite dimensional
subspace of $\dot{H}^{1}$, \bea \nonumber V_m = \hbox{Span} \{ e_k =
e^{ik.x},\; k=(k_1,k_2); 0 < \max(|k_1|, |k_2|) \leq m \}, \eea
endowed with the same inner product and the same norm as those of
$\dot{H}^1.$\\
Accordingly, we denote by $P_m$, the orthogonal projection onto
$V_m$, defined by: \bea\non P_m : \dot{L}^2 &\rightarrow & V_m \\
\non u &\mapsto & P_m(u) = \sum_{\max(|k_1|, |k_2|)\leq m}
\hat{u}(k) e_k ,\eea which commutes with the fractional Laplace
operator.

\fin \\

\section{Proof of Theorem \ref{TH1}}
We shall split the work into three steps : existence, uniqueness and
regularity for solutions of \Ref{imp1}. To begin with, we prove the
first step.

\subsection{Existence}

\noindent We take the inner product of \Ref{imp1} with
$\theta^{n+1}$ in $\dot{L}^{2}$. Using periodic boundary conditions
together with the fact that $u^{n+1}$ is divergence free, this leads
to: \bea \label{lin2}
\parallel\theta^{n+1}\parallel_2^2 \;  + \nu \tau \; \parallel
\Lambda^{\alpha}\theta^{n+1}\parallel_2^2 =
 (\tau f + \theta^{n},\theta^{n+1}).
  \eea
Thus, by Young's and Cauchy Schwartz's inequalities, we obtain: \bea
\label{lin3} \parallel\theta^{n+1}\parallel_2^2 \;  + 2 \nu \tau \;
\parallel\Lambda^{\alpha}\theta^{n+1}\parallel_2^2 \; \leq \;
 2 \tau^2 \; \parallel f\parallel_2^2 \; + \; 2 \parallel\theta^{n}\parallel_2^2.
  \eea
We infer from the a priori estimate \Ref{lin3}, that we shall look
for a weak solution $\theta^{n+1}$ that belongs to
$\dot{H}^{\alpha}$. On the other hand, by considering a variational
formulation of our problem, we remark that nonlinearity of the
variational form under consideration prevents us from resolving the
equation \Ref{imp1} by Lax-Milgram Lemma in $ \dot{H}^{\alpha}$.
Therefore, to contribute to the control of the nonlinearity, we
mimic the strategy in \cite{cordoba2004} for (QG) equation, so we
proceed to a variational regularization of \Ref{imp1}, which reads:
 \bea\label{imp2} & for\;\varepsilon>0,\;\;\;find
\;\;\theta_{\varepsilon}^{n+1}\in
\dot{H}^1\;\;such\;\;that,\;\;\forall v\in \dot{H}^1&\\
\non &\frac {1}{\tau}(\theta_{\varepsilon}^{n+1}-\theta^n ,v) + \nu
( \Lambda^{\alpha} \theta_{\varepsilon}^{n+1}, \Lambda^{\alpha}v) +
(\nabla.(u_{\varepsilon}^{n+1}\theta_{\varepsilon}^{n+1} ),v)&\\
\non &+ \varepsilon (\nabla\theta_{\varepsilon}^{n+1},\nabla v) =
(f,v).& \eea with $ u_{\varepsilon}^{n+1} = \; {\mathcal
R^{\bot}}(\theta_{\varepsilon}^{n+1}).$

\noindent  The problem \Ref{imp2} is nonlinear and its resolution is
based on the Faedo-Galerkin approximation method introduced in
Section 2.

\noindent To approach  $\theta_{\varepsilon}^{n+1}$ the solution of
\Ref{imp2}, we have to solve the following variational problem:
 \bea\label{imp3} &for\;\varepsilon>0,\;and\; m\in \mathbb{N}^{\ast},\;\;\;find
\;\;\theta_{\varepsilon,m}^{n+1}\in
V_m\;\;such\;\;that,\;\;\forall v\in V_m &\\
\nonumber &\frac {1}{\tau}(\theta_{\varepsilon,m}^{n+1}-P_m\theta^n
,v) + \nu  ( \Lambda^{\alpha} \theta_{\varepsilon , m}^{n+1},
\Lambda^{\alpha}v)  +
(\nabla.(u_{\varepsilon,m}^{n+1}\theta_{\varepsilon,m}^{n+1} ),v)&\\
\non &+ \varepsilon (\nabla\theta_{\varepsilon,m}^{n+1},\nabla v) =
(P_m f,v).& \eea We state and prove the following result:
\begin{prop}\label{prop3}
$ \forall m \in \N^{*}\;\;and\;\;\forall\varepsilon>0,\;$ there
exists $\theta_{\varepsilon,m}^{n+1} \in V_m$ a solution of
\Ref{imp3}. Moreover, we have the following a priori estimates:
\bea\label{1}
\parallel \theta_{\varepsilon,m}^{n+1}\parallel_2 \leq 2K_0,\\
\label{2}  \parallel\Lambda^{\alpha}\theta_{\varepsilon,m}^{n+1}\parallel_2 \leq \frac{K_0}{\sqrt{\nu\tau}},\\
\label{3} \parallel\nabla\theta_{\varepsilon,m}^{n+1}\parallel_2
\leq \frac{K_0}{\sqrt{\epsilon\tau}} \eea where $\; K_0 =
\sqrt{\tau^{2}\parallel f
\parallel_{2}^{2}+\parallel\theta^{n}\parallel_{2}^2 + 1}.\;$
\end{prop}
{\bf Proof : } in order to prove the existence of such solution, we
need the technical Brouwer's Lemma \ref{BrLem}. For that purpose,
consider here $X = V_m$, and $F : V_m \rightarrow V_m $ be defined
by
 \bea
\label{imp4} F(\theta_m)= \theta_{m} + \nu\tau \Lambda^{2\alpha}
\theta_{m} + \tau P_m\nabla.(u_m\theta_m )  - \tau \varepsilon
\Delta \theta_m - \tau P_m f - P_m \theta^n. \eea First of all we
shall verify the conditions of Lemma \ref{BrLem}, on $\;F$ defined
by \Ref{imp4}. Proving the continuity of $\;F$ is straightforward
from the continuity of the operators $ \Lambda^{\alpha},\; {\mathcal
R},\; \nabla $ and $P_m$ on $V_m.$ Next, taking the $\dot{L}^2$
inner product of $F(\theta^{n+1}_{\varepsilon,m})$ with
$\theta^{n+1}_{\varepsilon,m}$, we get by  Young inequality:
\begin{eqnarray}
\label{F1}
[F(\theta^{n+1}_{\varepsilon,m}),\theta^{n+1}_{\varepsilon,m}]&\geq
&
\frac{1}{2}\parallel\theta^{n+1}_{\varepsilon,m}\parallel^{2}_{2}+\nu\tau\parallel\Lambda^{\alpha}\theta^{n+1}_{\varepsilon,m}\parallel^{2}_{2}
+\varepsilon\tau\parallel\nabla \theta^{n+1}_{\varepsilon,m}\parallel^{2}_{2} - K_{0}^2\\
\label{F2}&\geq & \frac 12
\parallel\theta^{n+1}_{\varepsilon,m}\parallel^{2}_{2} - K_{0}^2.
\end{eqnarray}
\noindent Thus, thanks to \Ref{F2}, it becomes clear that for $\;
\theta^{n+1}_{\varepsilon,m} \in V_m $ such that\\ $ \;
\parallel\theta^{n+1}_{\varepsilon,m}\parallel_{2}  = 2 K_0,\; $ we have:
\begin{eqnarray*}
[F(\theta^{n+1}_{\varepsilon,m}),\theta^{n+1}_{\varepsilon,m}]& \geq
&
 K_{0}^2 > 0.
\end{eqnarray*}
Hence by Brouwer's Lemma \ref{BrLem}, we obtain the existence of
$\theta^{n+1}_{\varepsilon,m} \in V_m\;$ such that $\;
\parallel\theta^{n+1}_{\varepsilon,m} \parallel_2 \leq 2 K_0, $ and
$\;F(\theta^{n+1}_{\varepsilon,m}) = 0 .\;$ Moreover, \Ref{2} and
\Ref{3} follow immediately from \Ref{F1}.
\fin \\

\noindent Now we state and prove:
\begin{prop}\label{prop4}
For all $ \varepsilon > 0,$ \Ref{imp2} admits  a solution
$\theta_{\varepsilon}^{n+1} \in \dot{H}^{1}$ that satisfies \Ref{1},
\Ref{2} and \Ref{3}.
\end{prop}
{\bf Proof : } obviously, such a result is obtained by getting the
limit on m. At first sight, the estimates \Ref{1} and \Ref{3} infer
that $(\theta^{n+1}_{\varepsilon,m})_{m}$ is bounded in
$\dot{H}^{1}$ then it admits a subsequence still denoted by
$(\theta^{n+1}_{\varepsilon,m})_{m}$ which weakly converges to
$\theta_{\varepsilon}^{n+1}$ in $\dot{H}^{1}$, and  strongly in
$\dot{H}^{\alpha}$ and in $\dot{L}^4$, owing to the compact Sobolev
imbedding \bea \label{inj} \dot{H}^1\hookrightarrow \dot{H}^{\alpha}
\hookrightarrow \dot{L}^4.\eea Thus, we go back to \Ref{imp3} and we
let m goes towards the infinity. Using the continuity of the Riesz
operator on $\dot{L}^p$ spaces, we get
\begin{eqnarray}
u_{\varepsilon,m}^{n+1}&\rightarrow & u_{\varepsilon}^{n+1} =
{\mathcal R^{\bot}}(\theta_{\varepsilon}^{n+1} )\quad
in\;\dot{L}^{4}.
\end{eqnarray}
From \Ref{inj} we deduce that,
\begin{eqnarray}
\Lambda^{\alpha}\theta_{\varepsilon,m}^{n+1}&\rightarrow &
\Lambda^{\alpha}\theta_{\varepsilon}^{n+1}\;\; in\;\dot{L}^{2}.
\end{eqnarray}
The same above arguments yield:
 \begin{eqnarray*}<\nabla.(
u_{\varepsilon,m}^{n+1}\theta_{\varepsilon,m}^{n+1}),\eta>_{(\dot{H}^{-1},\dot{H}^{1})}&\rightarrow
&<\nabla .( u_{\varepsilon}^{n+1}
\theta_{\varepsilon}^{n+1}),\eta>_{(\dot{H}^{-1},\dot{H}^{1})}.\\
\end{eqnarray*}
Thus, we conclude that $\theta_{\varepsilon}^{n+1}$ is a solution of
\Ref{imp2} which belongs to $\dot{H}^{1}$, and the estimates
\Ref{1}, \Ref{2} and \Ref{3}, follow promptly by getting the limits
on m.

\fin\\

\noindent Now we are ready to state the existence of a solution to
\Ref{imp1}.
\begin{prop}\label{prop4'}
For all $n\geq 0$, there exists $\theta^{n+1}$ solution of
\Ref{imp1}, which belongs to $\dot{H}^{\alpha}$.
\end{prop}
{\bf Proof}: using \Ref{1}, \Ref{2} and \Ref{3} on
$\theta_{\varepsilon}^{n+1}$, there exists a subsequence still
denoted by $( \theta_{\varepsilon}^{n+1})_{\varepsilon>0}$ such that
$\sqrt{\varepsilon}\nabla\theta^{n+1}_{\varepsilon} \rightharpoonup
h\quad  \hbox{in} \;\dot{L}^{2}$, $
\theta^{n+1}_{\varepsilon}\rightharpoonup \theta^{n+1}\quad
\hbox{in} \;\dot{H}^{\alpha}$, and $\theta^{n+1}_{\varepsilon} \;
\rightarrow \theta^{n+1}\; \hbox{in} \; \;\dot{L}^{2}$ when
$\varepsilon \rightarrow 0$.

\noindent Hence, on the one hand  $\varepsilon\Delta
\theta_{\varepsilon}^{n+1}\rightarrow  0 \quad \hbox{in}\;
\dot{H}^{-1}$, and on the other hand
$u_{\varepsilon}^{n+1}\theta_{\varepsilon}^{n+1}\rightarrow u^{n+1}
\theta^{n+1}\quad \hbox{in}\;\dot{L}^{2}$, so we get $\nabla(
u_{\varepsilon}^{n+1}\theta_{\varepsilon}^{n+1} )\rightarrow \nabla
( u^{n+1}\theta^{n+1}) \quad \hbox{in} \;\dot{H}^{-1}$ when
$\varepsilon \rightarrow 0$. Therefore, the proposition is proved.

\fin \\

\subsection{Regularity}
Actually we move to state regularity results. We point out that the
regularity of $\theta_{\varepsilon}^{n+1}$, the solution of
\Ref{imp2}, echoes directly on the regularity of $\theta^{n+1}$.
Thus, to begin with, we prove the following proposition, which
states some regularity results for $\theta_{\varepsilon}^{n+1}$.

\begin{prop} \label{propreg}
$\forall\;\varepsilon\;> 0,\;\forall\;n\in \N$,
$\theta_{\varepsilon}^{n+1}$ the solution of \Ref{imp2} belongs to
$\dot{H}^2.\;$ Furthermore, there exists $\; C_n > 0\;$ independent
of $\varepsilon$  such that \bea \label{ez1}
\|\theta_{\varepsilon}^{n+1}\|_{\dot{H}^{2\alpha}}\leq C_n ,\eea for
all $\alpha > \frac 23.\; $
\end{prop}
\noindent {\bf Proof :} let $B$ be the operator defined by
$$B:= (I-\varepsilon\tau\Delta)^{-1},$$
which is a regularizing operator of order 2. Then \Ref{imp2} can be
rewritten equivalently as \bea \label{imp5}
\theta_{\varepsilon}^{n+1} + \tau\nu B\adelta
\theta_{\varepsilon}^{n+1} + \tau
B\nabla(u_{\varepsilon}^{n+1}\theta_{\varepsilon}^{n+1} ) = B(\tau f
+ \theta^n).\eea Now, since $f \in \dot{L}^2$, and referring to
Eq.\Ref{imp5}, then the maximal regularity of
$\theta_{\varepsilon}^{n+1} $ is $\dot{H}^2.$ Therefore, we are
going to prove that $\theta_{\varepsilon}^{n+1} \in \dot{H}^2 $ in
two steps. Firstly, we recall that the regularity $\dot{H}^{1}$ for
$\theta_{\varepsilon}^{n+1}$ is ensured due to results of
Proposition \ref{prop4}. Then we remark that thanks to the Sobolev
imbeddings, and the continuity of the Riesz operator, we have $\;
u_{\varepsilon}^{n+1} \theta_{\varepsilon}^{n+1} \in \dot{L}^p,\; $
for all $ p \in [ 1, \infty).$ Thus, \bea
B\nabla(u_{\varepsilon}^{n+1}\theta_{\varepsilon}^{n+1} ) \in
\dot{H}^{1,p}. \eea Moreover, \bea
 B\adelta \theta_{\varepsilon}^{n+1}\;\in \dot{H}^{3-2\alpha}.
\eea On the other hand, due to the Sobolev imbeddings, $\;
\dot{H}^{3-2\alpha} \hookrightarrow \dot{H}^{1,\frac{2}{2\alpha
-1}}, \;$ when $\;\alpha < 1,$ and $\; \dot{H}^2 \hookrightarrow
\dot{H}^{1,p},\quad \forall p > 2,\;$ we deduce that for $\alpha<1$
 $ \; \theta_{\varepsilon}^{n+1} \in
\dot{H}^{1,p_0}, \;$ which is an algebra for some $p_0 > 2.$
Now, by a bootstrap argument and using the Sobolev imbeddings, we deduce that $ \; \theta_{\varepsilon}^{n+1} \in \dot{H}^2 \;$\\

\noindent Actually, we move on toward  the estimate \Ref{ez1}. For
that purpose, let $\beta$ be a real that satisfies $0 <
\beta\leq\alpha,$ and that have to be fixed later. We are going to
prove that, \bea \label{ez01}
\parallel\theta_{\varepsilon}^{n+1}\parallel_{\dot{H}^{\beta+\alpha}}\leq C_n ,\eea
for all $\alpha > \frac 23 $  and $\; 0  < \; \beta \leq 3
\alpha -  2.\;$\\
 For that purpose, we take $v=\Lambda^{2\beta}\theta_{\varepsilon}^{n+1}$ in
\Ref{imp2}, then we get: \bea \label{lin10}
&(\theta_{\varepsilon}^{n+1},\Lambda^{2\beta}\theta_{\varepsilon}^{n+1})
+ \nu\tau (\Lambda^{2\alpha}
\theta_{\varepsilon}^{n+1},\Lambda^{2\beta}\theta_{\varepsilon}^{n+1})
+ \tau(\nabla.(u_{\varepsilon}^{n+1}\theta_{\varepsilon}^{n+1}
),\Lambda^{2\beta}\theta_{\varepsilon}^{n+1})&\\ \nonumber &-
\varepsilon\tau (\Lambda^2
\theta_{\varepsilon}^{n+1},\Lambda^{2\beta}\theta_{\varepsilon}^{n+1})
 = (\tau
f+\theta^n,\Lambda^{2\beta}\theta_{\varepsilon}^{n+1}),& \eea which
leads to: \bea \nonumber
&\parallel\Lambda^{\beta}\theta_{\varepsilon}^{n+1}\parallel_{2}^{2}+\nu\tau
\parallel\Lambda^{\alpha +\beta} \theta_{\varepsilon}^{n+1}\parallel_{2}^{2} +
\varepsilon\tau\parallel\Lambda^{1+\beta}\theta_{\varepsilon}^{n+1}\parallel_{2}^{2}&\\
\label{lin11}&\leq  |(\tau
f+\theta^n,\Lambda^{2\beta}\theta_{\varepsilon}^{n+1})| +
\tau|(\nabla.(u_{\varepsilon}^{n+1}\theta_{\varepsilon}^{n+1}
),\Lambda^{2\beta}\theta_{\varepsilon}^{n+1})|.&\eea Using Cauchy
Schwartz and Young inequalities, together with the embedding \bea
\label{lin12} \dot{H}^{\alpha+\beta}\hookrightarrow \dot{H}^{2\beta}
\;\;\;\;for \;\;\beta \leq \alpha, \eea we obtain, \bea
\label{lin13} |(\tau
f+\theta^n,\Lambda^{2\beta}\theta_{\varepsilon}^{n+1})| \leq
\frac{C}{\nu\tau} \parallel\tau f+\theta^n\parallel_{2}^{2} + \frac
{\nu\tau}{4}
\parallel\Lambda^{\beta+\alpha}\theta_{\varepsilon}^{n+1}\parallel_{2}^{2},
\eea and similarly,
 \bea\label{lin14}
\tau|(\nabla.(u_{\varepsilon}^{n+1}\theta_{\varepsilon}^{n+1}
),\Lambda^{2\beta}\theta_{\varepsilon}^{n+1})| \leq  \frac
{\tau}{\nu}
\parallel\Lambda^{1+\beta-\alpha}(u_{\varepsilon}^{n+1}\theta_{\varepsilon}^{n+1})\parallel_{2}^{2}+
\frac {\tau\nu}{4}
\parallel\Lambda^{\beta+\alpha}\theta_{\varepsilon}^{n+1}\parallel_{2}^{2}.
\eea

\noindent In order to estimate the product
$u_{\varepsilon}^{n+1}\theta_{\varepsilon}^{n+1}$ in
$\dot{H}^{1+\beta-\alpha ,\;2}$, we shall use the product estimation
\Ref{KPV}, for $s= 1+\beta-\alpha, \;\;q=\frac{2}{1-\alpha}\;$, and
$\;p= \frac{2}{\alpha}$. By these data, the Sobolev imbedding \bea
\label{lin16} \dot{H}^{\alpha}\hookrightarrow \dot{L}^q
\;\;and\;\;\;\dot{H}^{\alpha}\hookrightarrow \dot{H}^{s,p},\eea are
satisfied for $\beta\leq 3\alpha-2$.

\noindent Now we make use of the continuity of the Riesz operator on
$\dot{L}^p$ to get:
\begin{equation}
\parallel\Lambda^s (u_{\varepsilon}^{n+1}\theta_{\varepsilon}^{n+1})\parallel_2\leq
c\parallel \theta_{\varepsilon}^{n+1}\parallel_q\parallel\Lambda^s
\theta_{\varepsilon}^{n+1}\parallel_p.
\end{equation}

\noindent Thus, by \Ref{lin16}, \Ref{2}, \Ref{lin13} and \Ref{lin14}
we derive the wished bound of $\theta_{\varepsilon}^{n+1}$, on the
$\dot{H}^{\beta+\alpha}$ norm. \\
Consequently, we obtain \Ref{ez1} using \Ref{ez01}and a bootstrap
argument.
\fin \\
\begin{prop}
For $\alpha \in ]0,\; 1],\;$ $\theta^{n+1}$ the solution defined by
Proposition \ref{prop4} belongs to $\dot{H}^{\alpha}.$ Moreover, for
$\alpha \in ] \frac 23, 1[,\;$ this solution belongs to
$\dot{H}^{2\alpha}.\; $
\end{prop}
{\bf Proof :} since the estimates \Ref{1} and \Ref{ez1} are
independent of $\varepsilon,\; $ then by making $\varepsilon
\rightarrow 0,\; $ we obtain the desired result.
\fin\\

\subsection{Uniqueness}

%
\begin{prop}\label{prop5} Let $\theta^n\in \dot{L}^2$ and
$\alpha\in\;[ \frac 23,\;1[$. Then, for $\tau > 0$ small enough,
$\theta^{n+1}$ the solution of \Ref{imp1} is unique.
\end{prop}

\noindent {\bf  Proof :} let $\theta_1^n,\; \theta_2^n \in
\dot{L}^2$ and consider $\theta_1^{n+1},\; \theta_2^{n+1} \in
\dot{H}^{\alpha}$ the respective solutions according to \Ref{imp1}.

\noindent Now, set $\theta^n = \theta_{2}^{n} - \theta_1^n, $
$\theta^{n+1} = \theta_{2}^{n+1} - \theta_1^{n+1}, $ and $u^{n+1} =
\mathcal{R}^{\bot} \theta^{n+1}.$ Then, $\theta^{n+1}$ satisfies,
\bea\label{diffl2} \theta^{n+1} - \theta^n + \nu \tau
(-\Delta)^{\alpha}\theta^{n+1} + \tau \nabla(u^{n+1} \theta_2^{n+1}
+ u_1^{n+1}\theta^{n+1}) = 0. \eea Taking the inner product of
\Ref{diffl2} with $\theta^{n+1} $, we find: \bea \label{u0}
\parallel\theta^{n+1} \parallel_2^2 - \parallel\theta^{n} \parallel_2^2 + 2\nu\tau  \parallel\Lambda^{\alpha} \theta^{n+1} \parallel_2^2 \leq
\mid\underbrace{2\tau \; \int\; \nabla( u^{n+1} \theta_2^{n+1} )
\theta^{n+1}}_{I_{n+1}}\mid. \eea Now, by Cauchy Schwartz
inequality, we find \bea |I_{n+1}| \leq 2 \tau
\parallel\Lambda^{1-\beta}( u^{n+1} \theta_2^{n+1}) \parallel_2
\parallel\Lambda^{\beta}\theta^{n+1} \parallel_2, \eea for some $ \beta \in [ \frac
23,\; \alpha[.$ Then, using the product estimation \Ref{KPV}, with
$s=1-\beta$, $\frac 1q =1-\beta$ and $\frac 1p = \beta - \frac 12$,
supplemented by the compact Sobolev Imbedding
$\dot{H}^{\alpha}\hookrightarrow \dot{H}^{\beta}$, and the
continuity of the Riesz operator in $\dot{L}^p$ spaces yield: \bea
\label{u2} |I_{n+1}| \leq 2C
\tau\parallel\theta_2^{n+1}\parallel_{\dot{H}^{\alpha}}
\parallel\Lambda^{\beta} \theta^{n+1}
\parallel_2^2. \eea
Moreover, by interpolation we get: \bea \label{yb}
 \parallel\Lambda^{\beta} \theta^{n+1} \parallel_2^{2} \leq C  \parallel \theta^{n+1} \parallel_2^{2(1 -\frac{\beta}{\alpha})}  \parallel\Lambda^{\alpha} \theta^{n+1} \parallel_2^{2\frac{\beta}{\alpha}}
.\eea Now, using the Young Inequality and inserting \Ref{yb} in
\Ref{u2}, we deduce that, \bea \label{u3} |I_{n+1}| \leq  \frac 12
\parallel\theta^{n+1}\parallel_2^2 + C \tau^2 \parallel\theta_2^{n+1}\parallel_{\dot{H}^{\alpha}}^2
\parallel\Lambda^{\alpha} \theta^{n+1} \parallel_2^2. \eea Replacing \Ref{u3} in
\Ref{u0}, we deduce that there exists $\tilde{C}
> 0$ independent of $n$ such that \bea \label{u4} \parallel\theta^{n+1}
\parallel_2^2   + \tau [ 2\nu - \tilde{C}\tau\parallel\theta_2^{n+1}\parallel_{\dot{H}^{\alpha}}^2  ]  \parallel\Lambda^{\alpha}
\theta^{n+1} \parallel_2^2   \leq  2 \parallel\theta^{n}
\parallel_2^2. \eea Therefore, for $\tau > 0$ small enough such that \bea \tau
\parallel\theta_2^{n+1}\parallel_{\dot{H}^{\alpha}}^2 \leq \frac{2\nu}{\tilde{C}},  \eea we infer from
\Ref{u4} that, \bea \label{u5} \parallel\theta^{n+1} \parallel_2^2
\leq  2 \parallel\theta^{n} \parallel_2^2 \eea which yields the
uniqueness of $\theta^{n+1}$ the solution of \Ref{imp1}.

\fin \\
\section{Proof of Theorem \ref{TH2}}

We move to prove the existence and the regularity of the global
attractor. We recall that Theorem \ref{TH1} provides us a
semi-discrete dynamical system, $(\dot{L}^{2}, (S^n)_{ n \in \N})$
for $\alpha \geq \frac 23$, for $\tau$ small enough, that is given
by mean of the following map
 \bea
\nonumber S : &\dot{L}^{2} &\rightarrow \dot{L}^{2}\\
\nonumber &\theta^n & \mapsto S\theta^n = \theta^{n+1}, \eea where
$\theta^{n+1}$ is the unique solution of \Ref{imp1}, when $\alpha
\geq \frac 23$. Notice that following recursively Eq.\Ref{imp1}, and
starting from $\theta^0$, we define the operator $ S^n:
\;\dot{L}^{2}\rightarrow\;\dot{L}^{2} $ such that $ S^n \theta^0 =
\theta^n$.

\noindent To go ahead, it is well known that to describe the long
time behavior of solutions to the so defined dynamical system, we
shall concentrate on the dynamics of some absorbing sets for the
semi-group introduced above. We recall that general results
concerning the existence of global attractors are given in the book
of R. Temam \cite[Chapter 1]{temam} for both continuous and discrete
dynamical systems. To get the existence of the global attractor, we
have to fulfill the conditions of the following proposition:

\begin{prop}\label{pat}
Let $H$ be a Hilbert or complete metric space and let $ S \;: \; H
\rightarrow\; H $ be a continuous map, that satisfies the following
properties : \ben
\item  there exists a bounded absorbing set $\; {\mathcal B} \subset H, \;$
such that \bea \forall \theta^0  \in H,\; \exists
n_0(\theta_0),\quad \forall n \geq n_0(\theta_0),\quad S^n \theta_0
 \; \in \;{\mathcal{B}}
\eea
\item $S^n$ is uniformly compact for $n$ large enough. It means that
 for every bounded set $B \subset H,\;$ the set $S^n B$ is relatively compact in $\; H.$
\een Then, there exists an invariant compact set $\mathcal{A}
\subset H,\;$  that attracts all trajectories $S^n\theta^0,\;$ for
all $\;\theta^0 \in H.\;$ More precisely, \bea \non S^n(\mathcal{A})
= \mathcal{A},\qquad \hbox{and}\quad \hbox{dist}(S^n
\theta_0,\mathcal{A} ) \rightarrow 0,\quad \hbox{when}\quad n
\rightarrow \infty. \eea Hence, $\mathcal{A}= \omega (\mathcal B)$,
the $\omega-$limit set of ${\mathcal B}$, is the global attractor
for the semi-group $(S^n)_n$.
\fin \\

\end{prop}

\begin{rem}\label{rem1}
We notice that to verify the second condition proposed in
Proposition \ref{pat}, we can show that the set $S^n {\mathcal B}$
is bounded in a space compactly imbedded in $\; H.$
\fin \\
\end{rem}

\noindent To begin with, we prove the existence of some absorbing
sets.

\subsection{The $\dot{L}^2$ absorbing set}

\begin{prop}\label{prop6} Let $\tau \in ] 0, 1[,$ and
let $f\in \dot{L}^{2}.$ Consider $M$ a real that satisfies $M >
M_0,$ where $\;M_0$ is defined by \Ref{M0}. Then the set
\bea\label{ATT1} E=\{\theta\in \dot{L}^{2};\;\;
\parallel \theta \parallel_2\leq M\},\eea is an absorbing  set
positively invariant for $S$, that is, for all $\theta^0 \in
\dot{L}^{2}$ there exists $n_0 > 0$ such that \bea \label{ATT1a}
\forall n \geq n_0,\quad S^n \theta^0 \in E, \eea and
\bea\label{ATT1b} S(E) \subset E. \eea
\end{prop}

\noindent {\bf Proof:} taking the $\dot{L}^2$ inner product of
\Ref{imp1} with $\theta^{n+1}$ leads to:  \bea\label{ATT2}
\parallel \theta^{n+1}\parallel_{2}^{2}- \parallel
\theta^{n}\parallel_{2}^{2}+\parallel
\theta^{n+1}-\theta^{n}\parallel_{2}^{2}+ 2\nu\tau\parallel
\Lambda^{\alpha}\theta^{n+1}\parallel_{2}^{2}= 2\tau
(f,\theta^{n+1}). \eea

\noindent Now, thanks to Cauchy-Schwartz and Young inequalities, we
find, \bea \label{f100}
 \mid(\; f,\theta^{n+1}\; )\mid \leq \frac{\parallel f\parallel_2^2}{2\nu C_0} + \frac{C_0\nu}{2}  \parallel \theta^{n+1}\parallel_{2}^{2}.
\eea Hence, inserting \Ref{f100} in the right hand side of
\Ref{ATT2}, and using the fact that $\parallel
\theta^{n+1}-\theta^{n}\parallel_{2}^{2}$ is a positif term,  we
get: \bea\label{ATT2a} (1 + \nu\tau C_0)\parallel
\theta^{n+1}\parallel_{2}^{2}+ 2\nu\tau\parallel
\Lambda^{\alpha}\theta^{n+1}\parallel_{2}^{2} \leq \parallel
\theta^{n}\parallel_{2}^{2}+\tau \frac{\parallel f\parallel_2^2}{\nu
C_0}. \eea So using the Poincaré inequality \Ref{poincaré}, we
obtain from \Ref{ATT2a}, \bea\label{ATT2b} ( 1 + \nu\tau C_0)
\parallel \theta^{n+1}\parallel_{2}^{2} \leq
\parallel \theta^{n}\parallel_{2}^{2} +  \frac{\tau \parallel
f\parallel_2^2}{\nu C_0}. \eea
 Now, we set $r =
\frac{1}{1+\nu\tau C_0 }\;,$ and we rewrite \Ref{ATT2b} as follows:
\bea \label{ATT3} \parallel \theta^{n+1}\parallel_{2}^{2}&\leq &
r\parallel \theta^{n}\parallel_{2}^{2} + r \frac{\tau}{\nu C_0}
\parallel f\parallel_2^2  . \eea Then, by a simple induction, we get
recursively from \Ref{ATT3}: \bea\label{ATT4}
\parallel \theta^{n}\parallel_{2}^{2} &\leq & r^n \parallel
\theta^0\parallel_{2}^{2}+ (1-r^n ) \frac{(1 + \nu \tau C_0)
\parallel f\parallel_2^2}{\nu^2 C_0^2}  .\eea We point out that $r<1\;$, thus, setting \bea \label{M0}
M_0^2 = \frac{(1 + \nu  C_0) \parallel f\parallel_2^2}{\nu^2 C_0^2},
\eea we conclude that for $M > M_0,$ the set $E$ defined by
\Ref{ATT1} satisfies \Ref{ATT1a}.

\noindent Moreover $E$ satisfies \Ref{ATT1b}. Indeed, let
$\theta^{n}$ belongs to $E$, then $$ \parallel
\theta^{n}\parallel_{2}^{2}\leq M^2. $$ On the other hand, since
$M_{0}<M\;$, and $$\frac{\tau}{\nu C_0} \parallel
f\parallel_{2}^{2}\leq \nu\tau C_0 M_{0}^{2},$$ \Ref{ATT3} leads to:
\bea \label{ATT7}
\parallel \theta^{n+1}\parallel_{2}^{2}&\leq & r(1
+\nu\tau C_0)M^{2}\\
&=& M^{2}. \eea Thus we obtain \Ref{ATT1b}.
\fin\\

\subsection{The $\dot{L}^{p}$ bounded absorbing set}

\noindent Now, let us prove that for $p >2$, we get a uniform
boundedness on $\parallel\theta^{n+1}\parallel_p$. Moreover, we show
that there exists an absorbing ball for $\theta^{n+1}$ in the
$\dot{L}^p$ space.

\begin{prop}\label{propLP}
Let  $ 2 < p \leq \frac{2}{1-\alpha},\;  $ and $ f \in \dot{L}^p.$
We have for all $n \geq 1,$ \bea
\label{eqLp}\parallel\theta^{n+1}\parallel_p \leq \frac{1}{1 + \frac
 {2}{p}\nu\tau C_0} ( \parallel\theta^n\parallel_p
+ \tau\parallel f\parallel_p). \eea Moreover, the set
\bea\label{setpalpha} {\mathcal G} =\{\theta\in
\dot{L}^{p_{\alpha}}\;;\;
\parallel\theta\parallel_{p_{\alpha}} \leq M\},\eea is a bounded absorbing set for $S$, where $p_{\alpha} = \frac{2}{1-\alpha}$ and $ M > M_1$ where
\bea \label{M1} M_1 \;=\; \frac{1}{(1-\alpha)\nu C_0}\parallel
f\parallel_{p_{\alpha}}. \eea
 That
is, for all  $\theta^0,\;f\; \in \dot{L}^{p_{\alpha}},$ there exists
$n_1
> 0$ such that, $\forall n \geq n_1$ \bea S^n \theta^0\in {\mathcal
G} .\eea
\end{prop}

\noindent{\bf Proof:} suppose that $ p > 2$, then we proceed as in
\cite[page 172]{Ju2}. Our aim is to show for a given (fixed)
$\theta^0$, that $\parallel\theta^{n+1}\parallel_p$ is also
uniformly bounded for $ t
> 0 $.

\noindent  By taking  the inner product of \Ref{imp1} with $ p
|\theta^{n+1}|^{p-2} \theta^{n+1}$ in $\dot{L}^2$, we get: \bea
\nonumber p\parallel\theta^{n+1}
\parallel_{p}^{p} - p\int_{\Omega} \theta^n |\theta^{n+1}|^{p-2}
\theta^{n+1} +
\nu \tau p \int_{\Omega} \Lambda^{2\alpha}\theta^{n+1} |\theta^{n+1}|^{p-2} \theta^{n+1}  + \\
\nonumber
p\tau \int_{\Omega} u^{n+1}. \nabla\theta^{n+1}  |\theta^{n+1}|^{p-2} \theta^{n+1} = p\tau \int_{\Omega} f  |\theta^{n+1}|^{p-2} \theta^{n+1}.\\
\label{normelp} \eea

\noindent By an integration by parts and using the fact that
$\nabla.u^{n+1} = 0,$ we get \bea \label{null} p\int_{\Omega} u^{n}.
\nabla\theta^{n+1} |\theta^{n+1}|^{p-2} \theta^{n+1} = -
p\int_{\Omega} \nabla.u^{n+1} |\theta^{n+1}|^{p} = 0. \eea

\begin{rem}
In order to apply the improved positivity Lemma, we have to consider
$\theta_{\varepsilon,m}^{n+1}$ defined by \Ref{imp3}, instead of
$\theta^{n+1}$, since it satisfies the required assumptions of Lemma
\ref{lemlp}. Hence, all the computations are made formally on
$\theta^{n+1}$, however they are valid using the formulation
\Ref{imp3}, as it is observed earlier, and we conclude by taking the
limit on $\varepsilon$.
\end{rem}
\noindent Now, we shall use the improved positivity Lemma
\ref{lemlp}, and we focus particularly on Eq.\Ref{positivity}, to
get: \bea \label{positivityalpha} p \int_{\Omega}
|\theta^{n+1}|^{p-2} \theta^{n+1} \Lambda^{2\alpha} \theta^{n+1}
\geq 2 \int_{\Omega} (\Lambda^{\alpha} |\theta^{n+1}|^{\frac p2})
^2. \eea

\noindent On the other hand, thanks to the spectral properties of
the operator $\Lambda$, we have  \bea \label{pc} \int_{\Omega}
(\Lambda^{\alpha} |\theta^{n+1}|^{\frac p2})^2 \geq C_0
\parallel\theta^{n+1}\parallel_{p}^{p}. \eea
 Now gathering \Ref{null},
\Ref{positivityalpha} and \Ref{pc} in \Ref{normelp} we obtain after
using the H{\"o}lder inequality and simplifying by $
\parallel\theta^{n+1} \parallel_p^{p-1}  $ that,

 \bea ( p + \frac2p \nu
\tau C_0 ) \parallel\theta^{n+1} \parallel_p \leq
\parallel\theta^{n} \parallel_p +  \tau \parallel f \parallel_p. \eea
Thus, \bea \label{theequation}(1 + \frac2p \nu \tau C_0 )
\parallel\theta^{n+1} \parallel_p \leq
\parallel\theta^{n} \parallel_p + \tau\parallel f \parallel_p. \eea

\noindent We set, \bea K=\frac{1}{1+\frac {2}{p} \nu\tau C_0}. \eea
By a simple induction on \Ref{theequation} we infer that,
 \bea\nonumber
\parallel \theta^{n}\parallel_{p} &\leq & K^{n} \parallel
\theta^0\parallel_{p}+\tau\sum_{k=1}^{n} K^{k} \parallel f\parallel_p,\\
\label{Lp2} & \leq &  K^{n} \parallel \theta^0\parallel_{p}+
(1-K^{n})\frac{p}{2\nu C_0}\parallel f\parallel_{p}.\eea We set,
\bea \tilde{M} = \frac{p}{2\nu C_0}\parallel f\parallel_{p}. \eea
Then, for $n\geq n_1(\parallel \theta^0\parallel_{p})$, we'll get
the uniform boundedness of $
\parallel\theta^{n} \parallel_p $ independently of $\theta^0.\;$ Thus,
since $M > M_1\geq \tilde{M}$, then the set \bea \label{F} F =
\left\{ \theta \in \dot{L}^p, \;
\parallel \theta
\parallel_p \leq M \right\} \eea is absorbing and positively
invariant by $S.$
\fin\\
\subsection{The $\dot{H}^{\alpha}$ absorbing set}

\noindent Now, to fulfill the second condition of Proposition
\ref{pat}, we have to prove the following result:

\begin{prop}\label{propalpha} Let $f\in \dot{L}^{p_{\alpha}},$ and $N>0$ an integer, then we set $r=N\tau$. Consider
$M  >  M_2 $ where \bea \label{M2} M_{2} = (\frac{r}{\nu}\parallel
f\parallel_2^2 + \frac{a_2}{r}) \exp (\frac{rC}{1-\tau rC}), \eea
and $a_2$ is given by \Ref{a2}, then the set \bea\label{ATTalpha}
{\mathcal{B}} \; = \; \{\theta\in   E ;\;\;
\parallel \Lambda^{\alpha}\theta \parallel_2\leq M\},\eea is a bounded absorbing  set
 for $S$.
\end{prop}
{\bf Proof :} taking the $\dot{L}^2$ inner product of \Ref{imp1}
with $\Lambda^{2\alpha}\theta^{n+1}$ leads to: \bea \nonumber&
\frac{1}{2\tau} [
\parallel\Lambda^{\alpha}\theta^{n+1} \parallel_2^2 - \parallel\Lambda^{\alpha}\theta^{n}  \parallel_2^2 +
\parallel\Lambda^{\alpha}(\theta^{n+1} - \theta^{n} )\parallel_2^2   ] +
\nu
\parallel\Lambda^{2\alpha} \theta^{n+1} \parallel_2^2&\\ \label{Hs}&= \int_{\Omega} f
\Lambda^{2\alpha} \theta^{n+1}+\int_{\Omega} \nabla(u^{n+1}
\theta^{n+1}) \Lambda^{2\alpha} \theta^{n+1}. &\eea At first, We
estimate the first term in the right hand side of \Ref{Hs}, using
Cauchy Schwartz and Young Inequalities. Thus we obtain:
\bea\label{Hs'}\mid \int_{\Omega} f \Lambda^{2\alpha}
\theta^{n+1}\mid \leq \frac{\parallel f\parallel_2^2}{\nu} +
\frac{\nu}{4}\parallel \Lambda^{2\alpha}
\theta^{n+1}\parallel_2^2.\eea Secondly, we shall estimate the
nonlinear part of \Ref{Hs}. Namely, For some $0<\beta \leq \alpha$
we have, \bea\label{Hs3} \mid \int_{\Omega} \nabla(u^{n+1}
\theta^{n+1}) \Lambda^{2\alpha} \theta^{n+1}\mid \leq C
\parallel\Lambda^{1 + \alpha - \beta} (u^{n+1} \theta^{n+1} )\parallel_2 \parallel\Lambda^{ \alpha
+\beta}  \theta^{n+1} \parallel_2. \eea Actually, we take use again
of the product estimate \Ref{KPV}, so we get for $\frac 1p + \frac
1q = \frac 12$ \bea
\parallel\Lambda^{1 + \alpha - \beta} (u^{n+1} \theta^{n+1} )\parallel_2  \leq  C [
\parallel\Lambda^{1 + \alpha - \beta} u^{n+1} \parallel_p \parallel \theta^{n+1} \parallel_q +
\parallel\Lambda^{1 + \alpha - \beta} \theta^{n+1} \parallel_p \parallel u^{n+1} \parallel_q ]. \eea
Using the continuity of the Riesz operator on $L^p$ spaces for $ 1 <
p < \infty$, and the sobolev imbedding \bea \dot{H}^{\alpha+\beta}
\subset \dot{H}^{1 + \alpha -\beta, p}, \eea for $\beta = 1 - \frac
1p,$ and for $\frac 1q = \beta - \frac 12,$ we deduce that:
\bea\label{Hs43}
\parallel\Lambda^{1 + \alpha - \beta} (u^{n+1} \theta^{n+1}
)\parallel_2  \leq C
\parallel \theta^{n+1}\parallel_q
\parallel\Lambda^{ \alpha +\beta} \theta^{n+1}
\parallel_2. \eea

\noindent Now, using the $\dot{L}^q$ uniform boundedness in
Proposition \ref{propLP}, having $2<q\leq \frac {2}{1-\alpha}$, we
deduce that \bea | \int_{\Omega} \nabla(u^{n+1} \theta^{n+1})
\Lambda^{2\alpha} \theta^{n+1}  | \leq C
\parallel\Lambda^{ \alpha +\beta} \theta^{n+1}
\parallel_2^2. \eea Having $0<\beta\leq \alpha$, an easy interpolation yields:\bea
\parallel\Lambda^{ \alpha +\beta}  \theta^{n+1} \parallel_2 \leq
\parallel\Lambda^{ \alpha}\theta^{n+1} \parallel_2^{1-\frac{\beta}{\alpha}}
\parallel\Lambda^{ 2\alpha}  \theta^{n+1}
\parallel_2^{\frac{\beta}{\alpha}},
\eea such that, by Young inequality we obtain \bea \label{Hs4}|
\int_{\Omega} \nabla(u^{n+1} \theta^{n+1}) \Lambda^{2\alpha}
\theta^{n+1}  | \leq C
\parallel\Lambda^{ \alpha } \theta^{n+1} \parallel_2^2 + \frac{\nu}{4}
\parallel\Lambda^{ 2\alpha} \theta^{n+1} \parallel_2^2. \eea
Then, we get by inserting \Ref{Hs4} and \Ref{Hs'} in \Ref{Hs}: \bea
\nonumber \frac{1}{\tau} [
\parallel\Lambda^{ \alpha }\theta^{n+1}  \parallel_2^2 &-&
\parallel\Lambda^{ \alpha }\theta^{n}
\parallel_2^2 ]\\\label{Hs2} &+&  \nu
\parallel\Lambda^{2\alpha} \theta^{n+1}  \parallel_{2}^2  \leq C
\parallel\Lambda^{ \alpha } \theta^{n+1} \parallel_2^2 + C_1 . \eea
By Proposition \ref{prop6}, we get \bea \tau \sum_{ n=n_0}^{N+n_0}
\parallel\Lambda^{\alpha} \theta^{n+1} \parallel_2^2 \leq a_2, \eea
where \bea\label{a2}a_2= \frac{r}{\nu C_0}\parallel f\parallel_2^2 +
r\nu C_0 M^2. \eea On the other hand, a little care to Eq. \Ref{Hs2}
gives: \bea\label{Hs32} \frac{1}{\tau} [
\parallel\Lambda^{ \alpha }\theta^{n+1}  \parallel_2^2 &-&
\parallel\Lambda^{ \alpha }\theta^{n}
\parallel_2^2 ] \leq \frac{C}{1-\tau C}
\parallel\Lambda^{ \alpha } \theta^{n} \parallel_2^2 + \frac{C_1}{1-\tau C}.
\eea We emphasize that since $0<\tau=O(1)$, we can ensure that $
1-\tau C>0$, and that owing to the inequalities \Ref{Hs4},
\Ref{Hs43}, \Ref{setpalpha} and the definition of $M_1$, then $C=
C(C_0, \nu).M_1.$

\noindent Now thanks to the uniform Gronwall Lemma \ref{DUGL}, we
get the uniform boundedness of $
\parallel\Lambda^{\alpha} \theta^{n}
\parallel_{2}^{2},$ \bea\label{alpha} \parallel\Lambda^{\alpha} \theta^{n+1}
\parallel_{2}^{2}\leq (\frac{r}{\nu}\parallel f\parallel_2^2 +
\frac{a_2}{r}) \exp (\frac{rC}{1-\tau rC}), \qquad \forall n\geq n_0
+ N, \eea and hence we get the existence of ${\mathcal{B}}$.

\fin \\

\noindent

\subsection{Existence and regularity of the global attractor}

Before proceeding further to apply Proposition \ref{pat}, let us
reorder the previous results to depict the convenient phase space
$H$, that will allows us to fulfill the Proposition's conditions.\\
We fix $\; M > \max\left( M_0, M_1, M_2 \right)\;$ where $M_0, M_1 $
and  $M_2$ are defined by \Ref{M0}, \Ref{M1} and \Ref{M2} and
consider the set \bea H =  \left\{\theta\in \dot{L}^{p_\alpha};\;
\parallel \theta
\parallel_{p_{\alpha}}  \leq M  \right\}.
\eea Then, by the definition of $M$, and  owing to Proposition
\ref{prop6},
 Proposition \ref{propLP}, and Proposition \ref{propalpha}, there exists  $n\geq max(n_0 + N,\;n_1)$
such that if $\theta^n \in H$ then $\theta^{n+1}= S\theta^n \;\in
H.$ Hence, we have $$S: H \rightarrow H$$ is well defined. We define
$(H,\;d)$ as a complete metric space endowed with the metric $d$
defined by the $\dot{L}^2$ norm. It remains to prove the continuity
of $S$ on $(H,\;d). $ Therefore, we state the following Lemma:
\begin{prop}\label{propcont}
S is a continuous map from H to H for $\alpha> \frac 23$.
\end{prop}

\noindent {\bf Proof:} let $\theta_1^n, \;\theta_2^n \in H$ such
that $\theta_1^{n+1}= S \theta_1^n$ and $\theta_2^{n+1}= S
\theta_2^n.$\\ We set $\theta^{n+1}=\theta_2^{n+1}-\theta_1^{n+1}$
and $u^{n+1}= \mathcal{R}^{\bot}\theta^{n+1}$. Then $\theta^{n+1}$
satisfies: \bea\label{cont1} \theta^{n+1} + \tau
B\nabla(u^{n+1}\theta_2^{n+1} + u_1^{n+1} \theta^{n+1})= B \theta^n,
\eea where B is the linear operator defined by
$B:=(I+\nu\tau\adelta)^{-1}$, moreover,  we recall that for $\tau$
small enough, we have
\begin{equation}\label{ope}
\forall s_1 < s_2,~\|(1+ \nu \tau
\adelta)^{-1}\|_{\mathcal{L}(H^{s_1}, H^{s_2})} \leq
\frac{c}{\tau^{\frac{s_2-s_1}{2\alpha}}},
\end{equation}
particularly, we have: \bea\label{ope1}\parallel
B\parallel_{\mathcal{L}(\dot{L}^2,\dot{L}^2)}&\leq& 1,\\
\parallel
B\parallel_{\mathcal{L}(\dot{H}^t,\dot{H}^1)}&\leq&
\frac{C}{\tau^{\frac{1-t}{2\alpha}}},\eea for all $t<0$ to be fixed
later.

\noindent Thus we check the following estimations: \bea \nonumber
\parallel \theta^{n+1}\parallel_2&\leq & \parallel
B\parallel_{\mathcal{L} (\dot{L}^2,\dot{L}^2)}\parallel
\theta^n\parallel_2 + \tau\parallel
B(u^{n+1}\theta_2^{n+1})\parallel_{\dot{H}^1}\\\nonumber &+&
\tau\parallel
B(u_1^{n+1}\theta^{n+1})\parallel_{\dot{H}^1}\\
\nonumber &\leq &\parallel \theta^n\parallel_2 +
C\tau^{1-\frac{1-t}{2\alpha}}[\parallel
u^{n+1}\theta_2^{n+1}\parallel_{\dot{H}^t}+ \parallel
u_1^{n+1}\theta^{n+1}\parallel_{\dot{H}^t}].\eea Now, choosing $q=
\frac 32<2$ and $t= \frac {q-2}{q}=-\frac 13$ such that
$\dot{L}^q\hookrightarrow \dot{H}^t$,
 and using the continuity of
the Riesz operator and Hölder inequality  we obtain: \bea
\nonumber
\parallel \theta^{n+1}\parallel_2&\leq & \parallel \theta^n\parallel_2 + C \tau^{1-\frac{1-t}{2\alpha}}
[  \parallel \theta_2^{n+1}\parallel_{\dot{L}^6}  + \parallel
\theta_1^{n+1}\parallel_{\dot{L}^6}]  \parallel
\theta^{n+1}\parallel_{\dot{L}^2},\eea since $\frac 1q= \frac 12 +
\frac 16.$\\ Then, using the uniform boundedness of
$\parallel\theta_2^{n+1}\parallel_{\dot{L}^6}$ and
$\parallel\theta_1^{n+1}\parallel_{\dot{L}^6}$ owing to Proposition
\ref{propLP}, the defined set $F$ given by \Ref{F}, and the fact
that $6 \in ]2,\; p_{\alpha}]$ for $\alpha\geq \frac 23$, we get:
\bea
\nonumber
(1-C M\tau^{1-\frac{1-t}{2\alpha}})
\parallel \theta^{n+1}\parallel_2 \leq  \parallel
\theta^n\parallel_2.\eea

\begin{rem} At this stage, we must emphasize that $1-\frac{1-t}{2\alpha}>0$ only for $\alpha > \frac
23$, and hence for $ \tau > 0$ small enough, we get \bea \nonumber C
M\tau^{1-\frac{1-t}{2\alpha}}  < 1. \eea
\end{rem}

\noindent This makes end to this proof.

\fin\\

\begin{rem}
 Notice that the result of Proposition \ref{propcont}, yields the
 uniqueness of $\theta^{n+1}$ solution of \Ref{imp1}, on $H$. Hence,
 this gives rise to the dynamical system $(H,\;(S^n)_{n\geq 0})$.
\end{rem}

%
%

%
\noindent On the other hand, $S^n$ satisfies:

\begin{prop}\label{prop8}
 $(S^n)_{n\in \N}\;$ is uniformly compact in $H.$
\end{prop}
{\bf Proof:} since $\mathcal {B}$ is a compact subset in $H$ then
Proposition \ref{propalpha} achieves the proof.
\\
\fin \\
\noindent Therefore, the following main result is proved.
\begin{prop}\label{prop9}
For $\alpha> \frac 23$, the dynamical system $(H, (S^n)_{n\in \N}),$
admits a global attractor $\mathcal{A}$, which is included in
$\dot{H}^{2\alpha}$.
\end{prop}
{\bf Proof:} the assumptions of Proposition \ref{pat} are satisfied
thanks to the Propositions \ref{propalpha}, \ref{propcont}, and
\ref{prop8}.  Thus  there exists an invariant compact set $
\mathcal{A}$  included in $H$, which is the global attractor \bea
\label{lasta} \mathcal{A}= \omega ({\mathcal B}) =
\displaystyle{\bigcap_{n}  \overline{\bigcup_{k\geq n}S^k {\mathcal
B} }^{H}}, \eea
 the $\omega$-limit set of ${\mathcal B}$, and where the closure in \Ref{lasta} is taken with respect to the $\dot{L}^2$ metric \\

\noindent Regularity of this attractor have to be ensured from the
invariance property $S \mathcal{A}= \mathcal{A}$, and regularity
results of subsection 3.2.

\fin \\

\noindent Furthermore, we state and prove the following result:
\begin{prop}
$\mathcal{A}$ is a compact set in $\dot{H}^{\alpha}$.
\end{prop}
\noindent {\bf Proof:} to prove the compactness of the attractor, we
rely on the J. Ball argument \cite{jball}. We proceed as follows:
let $\; ( \theta_j^{n+1})_j\; $ a sequence of points of $\;
\mathcal{A} \subset \dot{H}^{\alpha}\;$ and
$u_j^{n+1}=\mathcal{R}^{\bot}(\theta_j^{n+1}).$ Now, consider the
sequence $\; ( \theta_j^{n})_j\; $ such that
$$
\theta_j^{n} = S^n \theta_j,
$$
\noindent or equivalently, \bea \label{impj}
\frac{\theta_j^{n+1}-\theta_j^n}{\tau} + \nu (-\Delta)^{\alpha}
\theta_j^{n+1} +\nabla.( u_j^{n+1}\theta_j^{n+1}) = f. \eea
Where $\theta^0_j = \theta_j.$ \\

\noindent At this stage we consider $n_0$ such that $\forall n \geq
n_0,\; \theta_j^n \in {\mathcal B}_{\alpha}.$

\noindent We emphasize that referring to the previous results and
subsections, there exist subsequences still denoted by
$(\theta_{j}^{n})_{j}$ and $(\theta_{j}^{n+1})_{j}$, such that \bea
\label{weak}
\theta_{j}^{n}\rightarrow \theta^n \quad \hbox{and}\quad \theta_{j}^{n+1}\rightarrow \theta^{n+1} \quad \hbox{weakly in }\quad \dot{H}^{\alpha},\\
\theta_{j}^{n}\rightarrow \theta^n \quad \hbox{and}\quad
\theta_{j}^{n+1}\rightarrow \theta^{n+1} \quad \hbox{strongly in
}\quad \dot{L}^{p}, \eea for any $ p \in  [ 2, \frac {2}{1-\alpha}]
.$ Hence, the limits $\theta^n$ and $\theta^{n+1}$ satisfy: \bea
\label{impj1} \frac{\theta^{n+1}-\theta^n}{\tau} + \nu
(-\Delta)^{\alpha} \theta^{n+1} +\nabla.( u^{n+1}\theta^{n+1}) = f.
\eea

\noindent We aim to prove that the convergence holds strongly in
$\dot{H}^{\alpha}$.\\


\noindent Let $w_j^n = \theta_j^n\;-\;\theta^n $ and $r_j^n =
\mathcal{R}^{\bot} w_j^n.$ By subtracting \Ref{impj1} from
\Ref{impj} we find \bea \label{impwj} \frac{w_j^{n+1}- w_j^n}{\tau}
+ \nu (-\Delta)^{\alpha} w_j^{n+1} +\nabla.( r_j^{n+1}
\theta_j^{n+1} + u^{n+1} w_j^{n+1}) = 0. \eea Taking the
$\dot{L}^{2}$ inner product of \Ref{impwj} with $w_j^{n+1}$ we
obtain, \bea \frac{1}{\tau} (\parallel w_j^{n+1}  \parallel_2^2  -
\parallel w_j^{n}
\parallel_2^2 ) &+& 2\nu \parallel \Lambda^{\alpha} w_j^{n+1} \parallel_2^2\\
\non &\leq & 2 \underbrace{ \mid\int \nabla.( r_j^{n+1}
\theta_j^{n+1}) w_j^{n+1}\mid}_{I_j^{n+1}}. \eea Let $\frac 23 <
\beta < \alpha,$ such that we get \bea \label{sobinjbeta}
\dot{H}^{2\beta - 1} \subset \dot{H}^{1 - \beta}. \eea By the
Cauchy-Schwartz inequality, we have \bea I_j^{n+1}  \leq
\parallel\Lambda^{1-\beta}(r_j^{n+1} \theta_j^{n+1})
\parallel_2 \parallel\Lambda^{\beta} w_j^{n+1} \parallel_2 .\eea Then, thanks to the Sobolev imbedding
\Ref{sobinjbeta}, we get: \bea I_j^{n+1}  \leq
\parallel\Lambda^{2\beta-1} (r_j^{n+1} \theta_j^{n+1})\parallel_2 \parallel\Lambda^{\beta} w_j^{n+1}
\parallel_2. \eea Now using the pointwise product estimate in Sobolev spaces
$\dot{H}^s(\T),\;$ for\\ $ 0 < s < 1,$ we get \bea I_j^{n+1} \leq
[\parallel\Lambda^{\beta} r_j^{n+1} \parallel_2
\parallel\Lambda^{\beta} \theta_j^{n+1}\parallel_2 ]\parallel\Lambda^{\beta}
w_j^{n+1} \parallel_2, \eea which using the continuity of the Riesz
operator, and the uniform boundedness results of Proposition
\ref{propalpha}, yields: \bea I_j^{n+1} \leq C
\parallel\Lambda^{\beta} w_j^{n+1} \parallel_2^2. \eea A simple interpolation
of $\dot{H}^{\beta}$ in $\dot{L}^2$ and $\dot{H}^{\alpha}$, together
with the Young inequality lead to: \bea I_j^{n+1} \leq C  \parallel
w_j^{n+1}
\parallel_2^2 + \nu
\parallel\Lambda^{\alpha} w_j^{n+1} \parallel_2^2. \eea Now using
the fact that the sequences $(w_j^n)_j$ and $(w_j^{n+1})_j$ converge
to $0$ in $\dot{L}^2$ we obtain that $(w_j^{n+1})_j$ converges to
$0$ in $\dot{H}^{\alpha}.$ Hence we get the compactness of
$\mathcal{A}$ in $\dot{H}^{\alpha}.$
\fin \\

\end{document}